\documentclass[12pt]{article}
\usepackage{latexsym,amsmath,amssymb}

\newif\ifdviwin

\usepackage[latin1]{inputenc}
\usepackage[english]{babel}
\usepackage{indentfirst}
\usepackage[mathscr]{eucal}
\usepackage{amssymb,amsmath,amsfonts}
\usepackage{fancybox,fancyhdr}
\usepackage{graphicx}

\newif\ifdviwin

\dviwintrue

\def\cS{\mathcal{S}}

\def\cM{\mathcal{M}}

\let\8=\infty \let\0=\emptyset

\let\landa=\lambda
\let\alfa=\alpha

\let\parc=\partial

\def\ep{\varepsilon}
\def\sg{\sigma}

\def\landa{\lambda}

\def\flecha{\rightarrow}
\def\esiz{\langle}
\def\esde{\rangle}

\def\cte.{\mathop{\rm cte.}\nolimits}

\def\E{\mathbb{E}}
\def\M{\mathbb{M}}

\def\R{\mathbb{R}}

\def\C{\mathbb{C}}

\def\H{\mathbb{H}}
\def\S{\mathbb{S}}

\def\Ek{\mathbb{E}^{3} (\kappa,\tau)}

 \newtheorem{defi}{Definition}
 \newtheorem{teo}[defi]{Theorem}
 \newtheorem{pro}[defi]{Proposition}
 \newtheorem{cor}[defi]{Corollary}
 \newtheorem{lem}[defi]{Lemma}
 
 \newtheorem{remark}[defi]{Remark}

 \newenvironment{proof}{\rm \trivlist \item[\hskip \labelsep{\it
      Proof}:]}{\par\nopagebreak \hfill $\Box$ \endtrivlist}

\numberwithin{equation}{section}

\begin{document}
\mbox{}\vspace{0.4cm}\mbox{}

\begin{center}
\rule{15.2cm}{1.5pt}\vspace{0.5cm}

{\Large \bf The Bonnet problem for surfaces in
\\[0.3cm] homogeneous $3$-manifolds}\\ \vspace{0.5cm} {\large José A.
Gálvez$\mbox{}^a$, Antonio Martínez$\mbox{}^b$ and Pablo
Mira$\mbox{}^c$}\\ \vspace{0.3cm} \rule{15.2cm}{1.5pt}
\end{center}
  \vspace{1cm}
$\mbox{}^a$, $\mbox{}^b$ Departamento de Geometría y Topología, Universidad de Granada,
E-18071 Granada, Spain. \\ e-mail: jagalvez@ugr.es ; amartine@ugr.es
\vspace{0.2cm}

\noindent $\mbox{}^c$ Departamento de Matemática Aplicada y Estadística,
Universidad Politécnica de Cartagena, E-30203 Cartagena, Murcia, Spain. \\
e-mail: pablo.mira@upct.es \vspace{0.2cm}


\vspace{0.3cm}

 \begin{abstract}
We solve the Bonnet problem for surfaces in the homogeneous $3$-manifolds with a $4$-dimensional
isometry group. More specifically, we show that a simply connected real analytic surface in $\H^2\times \R$ or $\S^2\times \R$ is uniquely determined pointwise by its metric and its principal curvatures if and only if it is not a minimal or a properly helicoidal surface. In the remaining three types of homogeneous $3$-manifolds, we show that except for constant mean curvature surfaces and helicoidal surfaces, all simply connected real analytic surfaces are pointwise determined by their metric and principal curvatures.
  \end{abstract}

\vspace{0.5cm}

\section{Introduction}

The theory of surfaces in Riemannian homogeneous $3$-manifolds with a $4$-dimensional isometry group is currently experiencing a great development. The starting point of the renewed interest on these surfaces is the work \cite{AbRo} (see also \cite{Abr}), in which a holomorphic Hopf-type differential for CMC surfaces in these spaces was constructed. This fact suggested the possible existence of interesting results for surfaces in these homogeneous $3$-manifolds, and attracted the attention of many geometers. Some significant recent advances on this topic can be found in \cite{Abr,AbRo,ACT,AEG,BeTa,Dan,Dan2,FeMi1,FeMi2,HLR,HST,MoOn,SaE,SaTo}.

On the other hand, the consideration of a homogeneous $3$-manifold $\mathcal{M}^3$ with a $4$-dimensional isometry group  as an ambient space for a surface theory is among the most natural possibilities, apart from space forms. Indeed, they constitute highly symmetric $3$-manifolds, and are strongly related to the $3$-dimensional Thurston geometries. There are five distinct classes of them: the product spaces $\S^2 \times \R$ and $\H^2 \times \R$, the Heisenberg $3$-space ${\rm Nil_3}$, and manifolds with the isometry group of the Berger $3$-spheres and of the universal covering of ${\rm PSL} (2,\R)$. For more details, see \cite{Abr,Dan,FeMi2} and Section \ref{tres}.

The aim of the present work is to clarify the following question: \emph{how much geometric information of a surface in $\mathcal{M}^3$ is sufficient in order to determine the surface uniquely?} In \cite{Dan} B. Daniel proved that a surface in $\cM^3$ is uniquely determined by its first and second fundamental forms, and by the effect on the surface of a certain ambient Killing field. Moreover, these quantities satisfy an overdetermined system of PDEs, so it is reasonable to consider which data are necessary and which are superfluous in order to have the surface uniquely determined.

The specific problem that we will investigate is the following one : \emph{Are two surfaces in $\cM^3$ with the same first fundamental form and the same principal curvatures necessarily pointwise congruent?}

The answer to this question will be, roughly, the following one: almost every surface in $\cM^3$ is uniquely determined by its first fundamental form and its principal curvatures. However, there exist some exceptional surfaces for which this is not the case. Nevertheless, if we require real analyticity of the surfaces, they are necessarily constant mean curvature (CMC) surfaces, or are \emph{helicoidal}, i.e. invariant under a $1$-parameter isometry subgroup of the ambient space.

If we do not ask the surfaces to be real analytic, the solution to the problem is basically the same, except for the fact that now it is not forbidden to glue together in a smooth way open pieces of surfaces of different natures, what creates many unhandy situations. This is the only reason why we have restricted ourselves to the real analytic case.

We observe that, on the one hand, the existence of the above mentioned exceptional surfaces indicates that our hypothesis cannot be weakened, while on the other hand, the fact that the exceptional surfaces can be classified gives a satisfactory answer to the question we are discussing.

The previous problem is just the formulation in the homogeneous $3$-manifolds setting of the classical \emph{Bonnet problem} in $\R^3$. It asks whether the metric and the mean curvature (or equivalently, by Gauss' Theorema Egregium, the metric and the principal curvatures) are enough to determine a surface in $\R^3$ uniquely, and whether one can reach a classification of the surfaces which cannot be locally determined by their metric and principal curvatures, if any. This problem has been treated from a classical viewpoint, but also from a modern viewpoint, involving sophisticated techniques of integrable systems theory. A detailed report on both approaches, as well as its extension to space forms $\S^3$, $\H^3$, can be found in \cite{Bob,BoEt,Kam,KPP}.

The outline of the present work is the following one. In Section 2 we will give a detailed exposition of our main results. In Section 3 we will review, following
\cite{FeMi2,Dan}, the basic equations of immersed surfaces in homogeneous $3$-manifolds. In Section 4 we will provide examples of exceptional surfaces for which the Bonnet problem has a negative answer. Sections 5 and 6 will be devoted to show that in the real analytic category, these examples are the only possible exceptional surfaces for the Bonnet problem. Section 7 will contain some closing global remarks, and will indicate some natural open problems related to our results.

\section{The main results}

Let $\mathcal{M}^3$ denote a homogeneous $3$-manifold with a $4$-dimensional isometry group. We begin by addressing the question we wish to answer:

\vspace{0.3cm}

{\bf The Bonnet problem:} Let $\cS_1 ,\cS_2$ be two immersed oriented surfaces  in the homogeneous
$3$-manifold $\mathcal{M}^3$. Assume that there is a diffeomorphism $\Phi :\cS_1\flecha \cS_2$ which is a
local isometry between them, and so that the principal curvatures of $\cS_1$ at $p\in \cS_1$ agree
with the principal curvatures of $\cS_2$ at $\Phi(p)\in \cS_2$ for every $p\in \cS_1$. Is then
$\Phi$ the restriction to $\cS_1$ of an isometry $\Psi$ of $\mathcal{M}^3$?

\vspace{0.3cm}

This is the classical formulation of the Bonnet problem for surfaces in $\R^3$, but with the ambient $3$-space adapted to our setting. In order to formulate our results regarding the Bonnet problem in $\mathcal{M}^3$, we need first to distinguish between two kinds of congruencies that will appear. We shall say that two immersed surfaces $\psi,\psi^*:M^2\flecha \mathcal{M}^3$ are \emph{pointwise congruent} if there exists an isometry $\Psi$ of the ambient space $\mathcal{M}^3$ such that $\psi^* = \Psi \circ \psi$. We will say that $\psi$ and $\psi^*$ as above are \emph{globally congruent} (or simply \emph{congruent}) if there exists an isometry $\Psi$ of $\mathcal{M}^3$ and a diffeomorphism $\Gamma:M^2\flecha M^2$ such that $\psi^*\circ \Gamma = \Psi \circ \psi$.

Let us also introduce the following notion, again taken from the Bonnet problem in $\R^3$.

\begin{defi}
A \emph{Bonnet mate} of an immersed oriented surface $\cS_1$ in $\mathcal{M}^3$ is another oriented surface $\cS_2$ in $\mathcal{M}^3$ for which the answer to the question posed by the above Bonnet problem is negative. In that situation, the pair $\cS_1,\cS_2$ will be called a \emph{Bonnet pair} in $\mathcal{M}^3$.
\end{defi}
This definition immediately indicates that two Bonnet mates are never pointwise congruent. But nevertheless, they can be globally congruent, as we will see. This phenomenon also takes place for Bonnet pairs in $\R^3$, even though it is not explicitly remarked in general.

The solution to the Bonnet problem that we present here is composed by two different cases, depending on whether $\mathcal{M}^3$ is a product space $\M^2 (\kappa)\times \R$ or not. Here $\M^2(\kappa) =\S^2 (\kappa) $ or $\M^2(\kappa) =\H^2 (\kappa) $ depending on the sign of $\kappa$.

Let us expose first the case in which $\mathcal{M}^3$ is a product geometry $\M^2 (\kappa)\times \R$. It is then known that any simply connected minimal surface in $\mathcal{M}^3$ belongs to a continuous $1$-parameter family of minimal surfaces with the same induced metric and the same principal curvatures. This is called the \emph{associate family} of the minimal surface, in analogy with the Euclidean case. Its existence proves that for simply connected minimal surfaces in $\mathcal{M}^3$ (other than pieces of slices $\M^2 (\kappa)\times \{t_0\}$), the question posed by the Bonnet problem has a negative answer. In addition, we will show that there is another class of surfaces in $\mathcal{M}^3$ having a Bonnet mate: the simply connected pieces of \emph{properly helicoidal} surfaces in the product space $\M^2 (\kappa )\times \R$, i.e. surfaces invariant by a continuous $1$-parameter group of isometries of the ambient space not leaving the vertical axis pointwise fixed, and which are not right vertical cylinders over some curve in $\M^2 (\kappa)$. 


Once here, the solution to the Bonnet problem in $\M^2(\kappa)\times \R$ is the following one:
 \begin{teo}\label{main1}
Let $\mathcal{S}_1$, $\mathcal{S}_2$ be real analytic Bonnet mates in $\M^2(\kappa)\times \R$. Then both $\mathcal{S}_1$, $\mathcal{S}_2$ are either minimal surfaces or open pieces of properly helicoidal surfaces.

Conversely, any simply connected surface in $\M^2(\kappa)\times \R$ belonging to one of the following two families always has a Bonnet mate:
 \begin{enumerate}
 \item
Minimal surfaces in $\M^2(\kappa)\times \R$ that are not slices $\M^2(\kappa)\times \{t_0\}$.
 \item
Open pieces of properly helicoidal surfaces in $\M^2(\kappa)\times \R$.
 \end{enumerate}
 \end{teo}

We must remark here that any minimal surface in $\M^2(\kappa)\times \R$ and more generally, any CMC surface in $\mathcal{M}^3$, is real analytic.

\begin{remark}
The hypothesis of real analyticity in the theorem will be used only to ensure that the surfaces have an \emph{identity principle}, i.e. that two different surfaces cannot overlap over a common open set.

If we drop this assumption and work in the smooth category, the proof of the above theorem will still tell the following: \emph{if a surface $\mathcal{S}$ has a Bonnet mate, there exists a dense open set of $\mathcal{S}$ such that around any of its points the surface is minimal, or it is properly helicoidal, or it is pointwise congruent to its Bonnet mate.}
\end{remark}

Let us now describe our results in the case where $\mathcal{M}^3$ is \emph{not} a product geometry $\M^2 (\kappa)\times \R$. For that, we will need the following definition regarding the diffeomorphism $\Phi$ that appears in the formulation of the Bonnet problem: we will say that $\Phi$ \emph{preserves orientations} if $\{d\Phi (X), d\Phi (Y)\}$ is a positively oriented basis for $\cS_2$ whenever $\{X,Y\}$ is a positively oriented basis for $\cS_1$. Otherwise, we will say that $\Phi$ \emph{reverses orientations}. We shall also say that a Bonnet pair in $\mathcal{M}^3$ is \emph{positive} (resp. \emph{negative}) if its associated diffeomorphism $\Phi$ preserves (resp. reverses) orientations. Although this distinction is not necessary in $\R^3$, $\S^3$, $\H^3$ or $\M^2 (\kappa)\times \R$, it will be unavoidable in the present case. Let us also remark that if the Bonnet pair is minimal ($H=0$), it can always be assumed to be a positive Bonnet pair.

Once here, in order to describe the fundamental Bonnet pairs in $\mathcal{M}^3$ that will appear, let us consider the following geometric dualities. The first two of them are due to Daniel \cite{Dan}:
 \begin{enumerate}
 \item
Any simply connected CMC surface in $\mathcal{M}^3\neq \M^2 (\kappa)\times \R$ has an associate CMC surface in some $\M^2 (\kappa)\times \R$ with the same metric and conformal structure. This is a generalization of the usual Lawson-type correspondence for CMC surfaces in space forms. Two CMC surfaces related by this correspondence are called \emph{sister surfaces} (see Section \ref{tres} for more details).
 \item
Any simply connected surface in $\mathcal{M}^3\neq \M^2 (\kappa)\times \R$ with constant mean curvature $H\neq 0$ has an associated surface in $\mathcal{M}^3$ with the same metric and conformal structure, and with constant mean curvature $-H$. These surfaces, which are never pointwise congruent, are called \emph{twin immersions}.
 \item
Any non-minimal simply connected open piece of a helicoidal surface in $\mathcal{M}^3\neq \M^2 (\kappa)\times \R$ (i.e. a surface invariant under the action of a continuous $1$-parameter isometry subgroup of $\mathcal{M}^3$) has another piece of helicoidal surface associated to it, with the same metric and conformal structure, but with opposite mean curvature function. We call these surfaces, which are never pointwise congruent, \emph{helicoidal mates} in $\mathcal{M}^3$.
 \end{enumerate}
The solution to the Bonnet problem in the homogeneous $3$-manifolds that are not product geometries is the following one:

\begin{teo}\label{main2}
Let $\mathcal{M}^3$ denote a homogeneous $3$-manifold with a $4$-dimensional isometry group, so that $\mathcal{M}^3\neq \M^2 (\kappa )\times \R$.
 \begin{enumerate}
 \item
Let $\mathcal{S}_1$, $\cS_2$ denote two real analytic positive Bonnet mates in $\mathcal{M}^3$. Then $\cS_1$ and $\cS_2$ are open pieces of helicoidal CMC surfaces. As a matter of fact, they are the sister surfaces of two properly helicoidal CMC surfaces in a product space $\M^2 (\kappa)\times \R$. 

Conversely, any simply connected open piece of a helicoidal CMC surface in $\cM^3$ whose sister surface in $\M^2 (\kappa)\times \R$ is properly helicoidal has a positive Bonnet mate in $\cM^3$.
 \item
Let $\mathcal{S}_1$, $\cS_2$ denote two real analytic negative Bonnet mates in $\mathcal{M}^3$. Then both $\cS_1$, $\cS_2$ are either non-minimal CMC surfaces, or are open pieces of helicoidal surfaces in $\cM^3$. As a matter of fact, after a change of orientation in one of them, they are either CMC twin immersions, or are helicoidal mates in $\cM^3$

And conversely, any non-minimal simply connected surface in $\cM^3$ which is either a CMC surface or an open piece of a helicoidal surface, has a negative Bonnet mate in $\cM^3$.
 \end{enumerate}
\end{teo}

\begin{remark}
Again, the proof of Theorem \ref{main2} will actually show the following: \emph{any (non-analytic) surface $\mathcal{S}$ in $\mathcal{M}^3$ admitting a Bonnet mate, has a dense open set such that around any of its points the surface has constant mean curvature, or it is helicoidal, or it is pointwise congruent to its Bonnet mate.}
\end{remark}

\begin{remark}
Although two Bonnet mates are never pointwise congruent, we will see in
Section \ref{last} that they can sometimes be globally congruent.
\end{remark}

\section{Surface theory in homogeneous $3$-manifolds}\label{tres}

Despite there are five distinct classes of homogeneous $3$-manifolds with a $4$-dimensional
isometry group, it is possible to develop a unified treatment for all of them, by means of a pair
of real constants $(\kappa,\tau)$ verifying $\kappa \neq 4\tau^2$ (see \cite{Dan}).

Indeed, a homogeneous $3$-manifold with a $4$-dimensional isometry group can be seen as a fibration over the simply connected $2$-dimensional space form $\M^2 (\kappa)$ of constant curvature $\kappa$, so that: $(a)$ the fibers are geodesics, and $(b)$ translations along the fibers are isometries of the space.

Associated to this fibration we may define the \emph{vertical vector field} $\xi$ of $\Ek$ as the unit vector field associated to the translations
along the fibers. Moreover, as translations along the fibers are isometries of $\Ek$, the
vertical field $\xi$ is a Killing field.

The curvature $\kappa$ of the space form $\M^2 (\kappa)$ is the first real constant that will be used to classify these homogeneous $3$-manifolds. The other one will be the \emph{bundle curvature} $\tau$, which is the real constant such that
$ \overline{\nabla}_X \xi=\tau (X\times\xi)$ holds for any vector field $X$ on the manifold. Here $\overline{\nabla}$ is the Levi-Civita connection of the manifold and $\times$ denotes the cross product. In addition, the pair $(\kappa, \tau)$ satisfies $\kappa - 4\tau^2 \neq 0$.

With this, we have:
 \begin{enumerate}
 \item
If $\tau =0$ we get the product space $\M^2 (\kappa)\times \R$, that is, $\S^2 (\kappa)\times \R$ or $\H^2 (\kappa )\times \R$, depending on the sign of $\kappa$.
 \item
If $\kappa =0$ and $\tau \neq 0$, we get the Heisenberg $3$-space ${\rm Nil}_3$.
 \item
If $\kappa >0$ and $\tau \neq 0$, the homogeneous $3$-manifold has the isometry group of the Berger spheres $\S_{\rm ber}^3$.
 \item
If $\kappa <0$ and $\tau \neq 0$, we obtain the homogeneous $3$-manifolds whose isometry group is isomorphic to the one of the universal cover of ${\rm PSL} (2,\R)$.
 \end{enumerate}
We shall denote by $\Ek$ the homogeneous $3$-manifold described by the pair $(\kappa, \tau)$.

When $\tau =0$, the isometry group of the homogeneous manifold $\E (\kappa,0) \equiv \M^2 (\kappa)\times \R$ has four connected components. Indeed, an isometry of $\M^2 (\kappa)\times \R$ can preserve or reverse independently the orientation of the base and the fiber. Contrastingly, if $\tau \neq 0$, an isometry of $\Ek$ necessarily preserves or reverses simultaneously the orientations of both the base and the fiber. This will become important for our purposes later on.

For more details on homogeneous $3$-manifolds, the reader may consult for instance
\cite{Abr,BeTa,Dan,FeMi2} and references therein.

Now, let us consider an immersed surface $\psi:\Sigma\flecha \Ek$, and let us view $\Sigma$ as a Riemann surface with the conformal structure given by the induced metric $\esiz ,\esde$. Thus, if $z$ denotes a complex coordinate of $\Sigma$, we have $\esiz d\psi,d\psi\esde = \landa |dz|^2$ for some smooth positive function
$\landa$, that will be called the \emph{conformal factor} of $\psi$. Let us also define in terms of $z=s+i t$ the usual operators  $\parc_z = (\parc_s -i\parc_t )/2$ and $\parc_{\bar{z}} =(\parc_s +i \parc_t)/2$.

\begin{defi}\label{fundada} In the above setting, let $\eta$
be the unit normal map of $\psi$ in $\Ek$, and let $\xi$ denote the vertical unit Killing field of
$\E^3(\kappa,\tau)$. We will call the \emph{fundamental data} of $\psi$ to the globally defined uple
 \begin{equation}\label{fundamental}
 (\landa |dz|^2, u,H,p \, dz^2, A \, dz),
 \end{equation}
where

\vspace{0.2cm}

\def\arraystretch{1.1}\begin{tabular}{l}
\qquad$\lambda |dz|^2 := \esiz d\psi, d\psi\esde$ is the induced metric of $\psi$,\\
\qquad$u:\Sigma\flecha [-1,1]$  is the \emph{angle function} of $\psi$, i.e.
$u=\esiz \eta,\xi\esde,$\\
\qquad$H:\Sigma\flecha \R$  is the mean curvature of $\psi$,\\
\qquad$p\, dz^2 :=-\esiz \psi_z,\eta_z\esde \, dz^2 $ is the Hopf differential of $\psi$,\\
\qquad$A \, dz :=  \esiz \xi,\psi_z\esde \, dz $ is the $(1,0)$-part of the $1$-form $\esiz \xi, d\psi\esde$,
\end{tabular}

\vspace{0.2cm}
\noindent and $z$ is an arbitrary complex parameter for $\Sigma$.
\end{defi}

\begin{remark}
While $u,H$ are well defined functions on $\Sigma$, the quantities $\landa, p,A$ depend on the chosen parameter $z$ of $\Sigma$. Nevertheless, the metric $ds^2 =\landa |dz|^2$ as well as the Hopf differential $P= p dz^2$ and the tangent $1$-form $\mathcal{A} = A dz + \bar{A} d\bar{z}$ are well defined global objects on $\Sigma$.
\end{remark}

Let us also remark that when $\tau=0$ (i.e., $\Ek=\M^2(\kappa)\times\R$), we may write $\psi
=(N,h):\Sigma\flecha \M^2 (\kappa)\times \R$, where $h$ is the \emph{height function}. It then
turns out that $A=h_z$ in this setting.

The integrability equations for a surface in $\Ek$ are summarized in the following result, proved in \cite{FeMi2} (see also \cite{Dan}):

\begin{teo}\label{th:formulas}
The fundamental data of an immersed surface $\psi:\Sigma\flecha \Ek$ satisfy the following
integrability conditions with respect to an arbitrary complex parameter $z$ of the Riemann surface $\Sigma$:

\begin{equation}\label{lasces}
\left\{\def\arraystretch{1.3} \begin{array}{lccc} {\bf (C.1)} & p_{\bar{z}} & = & \displaystyle
\frac{\landa}{2} (H_z + u A (\kappa - 4\tau^2)). \\ {\bf (C.2)} & A_{\bar{z}} & = & \displaystyle
\frac{u \landa}{2} (H+i\tau) .\\ {\bf (C.3)} & u_{z} & = & - (H-i\tau) A -\displaystyle \frac{2
p}{\landa} \bar{A}.\\ {\bf (C.4)} & \displaystyle \frac{4 |A|^2}{\landa} & = & 1 - u^2 .
\end{array}\right.
\end{equation}

Conversely, let $\Sigma$ denote a simply connected Riemann surface, and consider on $\Sigma$ : (a) two smooth functions $H,u:\Sigma\flecha \R$ with $u^2 \leq 1$, (b) a Riemannian metric $ds^2$, (c) a complex $2$-form of $(2,0)$-type $P$, and (d) a real $1$-form $\mathcal{A}$. Assume furthermore that, if $z$ denotes an arbitrary complex parameter of $\Sigma$, and we write $ds^2 =\landa |dz|^2$, $P= p \, dz^2 $ and $\mathcal{A}= A dz + \bar{A}d\bar{z}$, then the set $(\landa,u,H,p,A)$ verify \eqref{lasces} for some real constants $\kappa,\tau$ with $\kappa - 4\tau^2 \neq 0$. Then there exists a surface $\psi:\Sigma\flecha \Ek$ with fundamental data \eqref{fundamental}. Moreover, this surface is unique up to isometries of $\Ek$ preserving the orientations of both the base and the fiber of $\Ek$.
%
%
%
%
\end{teo}

Equation ${\bf (C.1)}$ is nothing but the Codazzi equation of the ambient space. Even though these
four equations suffice in order to determine a surface in $\Ek$, we will also write down for later use two other
formulas that can be inferred from ${\bf (C.1)}$ to ${\bf (C.4)}$ (see \cite{Dan,FeMi2}). One is the \emph{Gauss equation}
 \begin{equation}\label{Gauss}
 K= {\rm det} (S) + \tau^2 + (\kappa -4\tau^2) u^2,
 \end{equation}
where $K$ is the Gaussian curvature of the metric $\landa |dz|^2$ and $S$ is the shape operator of
$\psi$. The other additional equation is
 \begin{equation*}
{\bf (C.0)} \hspace{1cm} A_z - \frac{\landa_z}{\landa} A = u p.
 \end{equation*}
It is also interesting to recall at this point the general relation
 \begin{equation}\label{dets}
 {\rm det} (S)= H^2 -\frac{4 |p|^2}{\landa^2},
 \end{equation}
that follows directly from the definition of $p,H$ and $\landa$.

\begin{remark}\label{isome}
Let $\psi:\Sigma\flecha \Ek$ be an immersed oriented surface with fundamental data \eqref{fundamental}, let $\Psi$ denote an isometry of $\Ek$, and consider the surface $\psi^* = \Psi\circ
\psi:\Sigma\flecha \Ek$, where we are taking for $\psi^*$ the orientation determined by the
complex structure of the Riemann surface $\Sigma$. Then, depending on $\Psi$ and the bundle
curvature $\tau$, the fundamental data of $\psi^*$ are the following ones:
 \begin{enumerate}
\item If $\Psi$ preserves the orientations of both the base and
the fiber of $\Ek$, the fundamental data of $\psi^*$ agree with the ones
of $\psi$. If $\Psi$ reverses both orientations, then $\psi^*$ has $(\landa |dz|^2,-u,H,p \, dz^2,-A \, dz)$ as its fundamental data.
 \item
If $\Psi$ preserves (resp. reverses) the orientation of the base and reverses (resp.
preserves) that of the fiber, then $\tau =0$ and the fundamental data of $\psi^*$ are $(\landa |dz|^2,u,-H,-p \, dz^2 ,-A\, dz)$ (resp.
$(\landa |dz|^2, -u,-H,-p\, dz^2,A \, dz)$).
 \end{enumerate}
\end{remark}

To close this section, we will expose a Lawson-type isometric correspondence between constant mean curvature (CMC) surfaces in different homogeneous $3$-manifolds. This correspondence is due to Daniel \cite{Dan}.

Let $\psi_1 :\Sigma\flecha \E^3 (\kappa_1,\tau_1)$ be a simply connected CMC surface in $\E^3 (\kappa_1,\tau_1)$ whose fundamental data are $(\landa |dz|^2,u,H_1,p_1\, dz^2 ,A_1\, dz)$, being $H_1$ constant. Consider now $\kappa_2,\tau_2 \in \R$ with $\kappa_1 - 4\tau_1^2 =\kappa_2 - 4\tau_2^2$, and take $H_2\in \R$ such that $H_2^2 + \tau_2^2 = H_1^2 + \tau_1^2$. So, there is some $\alfa\in \R$ such that $H_2 -i\tau_2 = e^{-i\alfa} (H_1- i\tau_1)$. It is then immediate to check that the set
\begin{equation}\label{fundadan}
(\landa |dz|^2, u, H_2, p_2= e^{i\alfa} p_1 \, dz^2, A_2=e^{i\alfa} A_1\, dz)
\end{equation}
 together with $\kappa_2, \tau_2$ verify conditions ${\bf (C.1)}$ to ${\bf (C.4)}$. So, by Theorem \ref{th:formulas} there exists a CMC surface $\psi_2 :\Sigma\flecha \E^3 (\kappa_2,\tau_2)$ that is locally isometric to $\psi_1$, and whose fundamental data are given by \eqref{fundadan}.

\section{Examples of Bonnet pairs in $\Ek$}\label{secmain}

In order to work with Bonnet pairs, we shall first define a common Riemann surface structure for any two Bonnet mates in $\Ek$. This will allow us to speak about common holomorphic objects of both surfaces.

Let $\cS_1,\cS_2$ be a Bonnet pair in $\Ek$, and let $\Phi:\cS_1\flecha \cS_2$ be the diffeomorphism in the formulation of the Bonnet problem.

If $\Phi$ preserves orientations,  $\cS_1$ and $\cS_2$ share the same conformal structure, as they have the same metric and orientation. Now, if $\Phi$ reverses orientations, and $\cS_2^*$ denotes the immersed surface $\cS_2$ with its opposite orientation, then $\cS_1$ and $\cS_2^*$ have the same conformal structure. However, they do not have the same principal curvatures at corresponding points anymore (except if $H=0$, see the Remark below). They have \emph{opposite} principal curvatures.

Nevertheless, let us assume in this second case that $\tau =0$. This allows us to consider $\Pi$ to be a vertical symmetry in $\M^2 (\kappa)\times \R$. Now define $\widetilde{\Phi} = \Pi \circ \Phi$ and $\widetilde{\cS_2}= \widetilde{\Phi} (\cS_2^*)$. As $\Phi$ is orientation reversing, we get that $\widetilde{\Phi}$ preserves orientations, and thus $\cS_1$ and $\widetilde{\cS_2}$ are in the condition of the first case.

In other words, if $\tau =0$ (and thus $\Ek \equiv \M^2 (\kappa)\times \R$), we can always assume in the formulation of the Bonnet problem that $\Phi$ preserves orientations. This is exactly what happens regarding the Bonnet problem in $\R^3$, $\S^3$ and $\H^3$, but it cannot be extended to the homogeneous spaces $\Ek$ with $\tau \neq 0$, as we will see.

\begin{remark}\label{lodeachecero}
If $\mathcal{S}_1$, $\mathcal{S}_2$ are two minimal surfaces in $\Ek$ that constitute a negative Bonnet pair, then $\mathcal{S}_1$ and $\mathcal{S}_2^*$ are a positive Bonnet pair with the same principal curvatures at every point. Thus, in the minimal case and up to a change of orientation, Bonnet pairs can always be assumed to be positive.
\end{remark}

Bearing the above facts in mind, we can conclude that the Bonnet problem for surfaces in $\Ek$ can be rephrased in the following way, which involves a common conformal structure for the two Bonnet mates:

 \begin{quote}
Let $\Sigma$ be a Riemann surface, and consider $\psi,\psi^*:\Sigma\flecha \Ek$ two conformal
immersions of $\Sigma$ into $\Ek$ with their respective orientations induced by the complex
structure of $\Sigma$. Assume that the induced metrics of $\psi$ and $\psi^*$ coincide, and that
$k_i =\ep k_i^*$ for $i=1,2$ and $\ep =\pm 1$ (with $\ep =1$ if $\tau =0$ or $H=0$), where $k_i$ (resp. $k_i^*$) are the principal
curvatures of $\psi$ (resp. $\psi^*$). Does it exist then an isometry $\Psi$ of $\Ek$ such that
$\Psi \circ \psi =\psi^*$?
 \end{quote}

\vspace{0.3cm}

Keeping the above formulation of the Bonnet problem in $\Ek$ in mind, our approach to a solution of
the problem will rely on the following lemma.

\begin{lem}\label{datos}
Let $\psi:\Sigma\flecha \Ek$ be a simply connected surface in $\Ek$ with fundamental data
\eqref{fundamental}. Assume that there exist a complex $2$-form $p^* \, dz^2$ and a complex $1$-form $A^* \, dz$ on $\Sigma$ with
\begin{equation}\label{noigual}
(p^* \, dz^2,A^*\, dz)\neq (\pm p \,dz^2,\pm A \, dz)
\end{equation}
and $\ep =\pm 1$ (with $\ep =1 $ if $\tau =0$) such that the set
 \begin{equation}\label{fundabon}
(\landa |dz|^2,u,\ep H,p^* \, dz^2,A^*\, dz)
 \end{equation}
verify conditions ${\bf (C.1)}$ to
${\bf (C.4)}$. Then $\psi$ has a Bonnet mate $\psi^*:\Sigma\flecha \Ek$.

On the other hand, if $\psi^*:\Sigma\flecha \Ek$ is a Bonnet mate of $\psi$, then its fundamental data are
of the form
 \begin{equation}\label{chulipiruli}
(\landa |dz|^2,u^*,\ep H,p^* \, dz^2,A^* \, dz)
 \end{equation}
for $\ep =\pm 1$ (with $\ep =1$ if $\tau =0$, up to a vertical symmetry), where $u^*:\Sigma\flecha [-1,1]$ verifies $u^2 = (u^* )^2$.
\end{lem}
\begin{proof}
Assume first of all the existence of a Bonnet mate $\psi^*:\Sigma\flecha \Ek$ of $\psi$, and let
$$(\landa^* |dz|^2,u^*,H^*,p^* \, dz^2,A^* \, dz)$$ denote its fundamental data. As the induced metric of both surfaces
agree, we have $\landa |dz|^2 =\landa^* |dz|^2$. By the above discussion, their principal curvatures agree up to
sign, so we get $H^*=\ep H$ for $\ep =\pm 1$. Now, by the Gauss equation \eqref{Gauss} we see that $u^2 = (u^*)^2$.

On the other hand, assume that we have two different sets of fundamental data of the form
\eqref{fundamental} and \eqref{fundabon} on the simply connected Riemann surface $\Sigma$,
so that $\ep =\pm 1$ and \eqref{noigual} holds. By Theorem \ref{th:formulas} and
\eqref{Gauss} they give rise to two immersions $\psi,\psi^*:\Sigma\flecha \Ek$ with the same
conformal structure, the same induced metric, and whose principal curvatures agree up to the $\ep$
sign. By making a change of orientation on $\psi^*$ if necessary (i.e. if $\ep =-1$) we can then
conclude that $\psi^*$ has the same metric and principal curvatures that $\psi$. So, we only need
to check that both immersions are non-congruent in $\Ek$, in order to have a Bonnet pair. But this
follows directly from the condition \eqref{noigual} and Remark \ref{isome}, what
finishes the proof.
\end{proof}

The remaining part of this section will be devoted to expose Bonnet pairs in homogeneous $3$-manifolds. It will be proved in Section 5 and Section 6 that, for real analytic surfaces, the
examples presented here are the only Bonnet pairs in $\Ek$.

\vspace{0.3cm}

{\bf Associate minimal surfaces in $\M^2 (\kappa)\times \R$:} Let $\psi :\Sigma \flecha \M^2 (\kappa)\times \R$ be a simply connected minimal surface, and write $ \psi =(N,h)$, where $h$ is the height function. Then $h:\Sigma\flecha \R$ is harmonic. Moreover, let $(\landa |dz|^2,u,0,p \, dz^2 ,A \, dz : =h_z \, dz)$ denote its fundamental data. It is then straightforward to check that for every $\theta \in [0,2\pi)$ the quantities $(\landa,u,0,e^{i\theta} p, e^{i\theta} A)$ verify conditions ${\bf (C.1)}$ to ${\bf (C.4)}$. Thus, by Theorem \ref{th:formulas} we get a $1$-parameter family $\psi_{\theta}$ of minimal immersions from $\Sigma$ into $\M^2 (\kappa)\times \R$ with the same induced metric and the same principal curvatures. Moreover, if $A\neq 0$ (i.e. the minimal surface is not a piece of a slice $\M^2 (\kappa)\times \{t_0\}$) and $\theta \neq 0,\pi$, the surfaces $\psi$ and $\psi_{\theta}$ are not pointwise congruent by Remark \ref{isome}. Thus, in general, any two elements of such an \emph{associated family} constitute a Bonnet pair in $\M^2 (\kappa)\times \R$.

\vspace{0.3cm}

{\bf Helicoidal surfaces in $\M^2 (\kappa)\times \R$:} Apart from the above associate minimal family, Bonnet pairs also appear when considering helicoidal surfaces in $\M^2 (\kappa)\times \R$, that is, surfaces that are invariant under a continuous $1$-parameter group of rigid motions of the ambient space. A general study of these surfaces can be consulted in \cite{MoOn,SaE,SaTo}. Our alternative approach here is based on the following characterization result:
 \begin{pro}\label{heli}
Let $\psi:\Sigma \flecha \Ek$ be a conformal immersion from the Riemann surface $\Sigma$, let $z$
be a conformal parameter for $\Sigma$, and suppose that $\psi$ has fundamental data
\eqref{fundamental}. If all the quantities $(\landa,u,H,p,A)$  depend only on $z+\bar{z}$, then
$\psi$ is an open piece of a helicoidal surface in $\Ek$.

Conversely, any helicoidal surface in $\Ek$ with fundamental data \eqref{fundamental} has around
any point a certain local conformal parameter $z$ for its induced complex structure such that
$(\landa,u,h,p,A)$ depend only on $z+\bar{z}$.
 \end{pro}
 \begin{proof}
Write $z=s+it$, and assume that $(\landa,u,H,p,A)$ depend only on $s$. In this way, we may consider the immersion $\psi$ to be defined on a vertical strip of $\C$. Then, for any
$t_0\in\mathbb{R}$ the map $(s,t)\mapsto (s,t+t_0)$ preserve the fundamental data of $\psi$, and
consequently $\psi$ has a continuous $1$-parameter group of self-congruences. In other words, for
every $t_0\in\mathbb{R}$ there is an isometry $\Psi_{t_0}$ of $\Ek$ satisfying
$\psi(s,t+t_0)=\Psi_{t_0} (\psi(s,t))$. So, $\{ \Psi_{t_0}\;:\;t_0\in\mathbb{R}\}$ is a continuous
$1$-parameter group of isometries of $\Ek$ i.e., it consists of helicoidal motions. Let us also
observe that if $\tau =0$, then $A=h_z$, and so the $1$-parameter isometry subgroup is made up by
rotations if and only if $h(s,t+t_0)=h(s,t)$, i.e. if and only if $A(s)\in \R$. Here a rotation means an isometry of $\Ek$ acting trivially on the fibers. In the same way,
$A(s)\in i\R$ if and only if the $1$-parameter isometry subgroup consists of vertical translations
and the surface is a piece of a right vertical cylinder.

Conversely, let $\psi:\Sigma\flecha \Ek$ be a surface invariant under a continuous $1$-parameter subgroup $\{\Psi_{v_2} :v_2\in \R\}$  of isometries of $\Ek$. If we choose $\alfa (v_1)$ a regular curve of the surface in $\Ek$ that is transversal to the orbits $\{\Psi_{v_2} (p) : v_2\in \R\}$ of the $1$-parameter subgroup, then we can parametrize the helicoidal surface $\psi$ locally as $\psi (v_1,v_2)= \Psi_{v_2} (\alfa (v_1))$ around any point that is not invariant under this isometry subgroup. Now, as $\psi (v_1,v_2+ \delta)$ is pointwise congruent to $\psi (v_1,v_2)$, we can conclude that the coefficients of the first and second fundamental forms of $\psi$ with respect to the coordinates $(v_1,v_2)$ do not depend on $v_2$, i.e. they depend on $v_1$ exclusively. The same holds regarding the function $u = \esiz \eta, \xi\esde$. Write $$ I= E(v_1) dv_1^2 + 2F(v_1) dv_1 dv_2 + G(v_1) dv_2^2$$ for the first fundamental form of $\psi$. Then we can consider the new local parameters $(s,t)$ given by $$s(v_1)= f (v_1)= \displaystyle \int^{v_1} \frac{\sqrt{EG-F^2}}{G}  dv_1, \hspace{0.7cm} t(v_1,v_2)= v_2 + \displaystyle\int^{v_1} \frac{F}{G} \ dv_1.$$ A direct computation lets us verify that  $I= G(v_1) (ds^2 +dt^2)$. And as $f$ is strictly increasing, we can write $ G(v_1)= G(f^{-1} (s)) $. At last, we obtain that
$$E=E(s) ,\hspace{1cm} F=0, \hspace{1cm} G=G(s) \hspace{1cm} E(s)=G(s),$$ and also that the coefficients of the second fundamental form of $\psi$ in the parameters $(s,t)$ depends only on $s$. Again, the same is true for $u$. All of this means that $z=s+it$ is a local conformal parameter of the surface, and all the coefficients of the fundamental data of $\psi$ only depend on $s=(z+\bar{z})/2$. This completes the proof.
 \end{proof}

By definition, a \emph{properly helicoidal} surface will be a helicoidal surface in $\M^2
(\kappa)\times \R$ whose generating $1$-parameter group of ambient isometries acts
non-trivially on both the vertical and horizontal factors. This excludes rotational
surfaces in $\M^2 (\kappa)\times \R$, as well as right vertical cylinders over some curve of $\M^2 (\kappa)\times \R$.

\begin{pro}\label{helibon}
Any simply connected open piece of a properly helicoidal surface in $\M^2 (\kappa)\times \R$ has a Bonnet mate.
\end{pro}
\begin{proof} Let $(\landa, H, u,p,A)$ denote the coefficients of the fundamental data of a properly helicoidal
surface $\psi :\Sigma\flecha \M^2 (\kappa)\times \R$, with respect to a local conformal parameter $z$ for it. By Proposition \ref{heli} we know that around each point of $\Sigma$ all these quantities can be assumed to depend exclusively on $s=(z+\bar{z})/2$. Now, when two such parameters $z_1,z_2$ coexist on an open set, the Cauchy-Riemann equations indicate that $z_2 = a z_1$ for some $a\in \R$ on this open set. So, as $\Sigma$ is simply connected, this implies the existence of a global conformal parameter $z$ on $\Sigma$ with the property that all the coefficients $(\landa, H, u,p,A)$ depend only on $s=(z+\bar{z})/2$. Now, a direct
computation shows that the quantities
 \begin{equation}\label{elidata}
 (\landa, u, H, \bar{p},\bar{A})
  \end{equation}
satisfy conditions ${\bf (C.1)} $ to ${\bf (C.4)}$ in $\M^2 (\kappa)\times \R$, and thus by Theorem
\ref{th:formulas} they give rise to an immersion $\psi^*:\Sigma\flecha \M^2 (\kappa)\times \R$
whose fundamental data are given by $(\landa |dz|^2, u, H, \bar{p}\, dz^2,\bar{A}\, dz)$. Thus,
$\psi$ and $\psi^*$ have the same metric and principal curvatures. Moreover, by the proof of
Proposition \ref{heli} we know that $A\notin \R$ and $A\notin i\R$ due to the fact that $\psi$ is
properly helicoidal. This indicates that $\psi$ and $\psi^*$ are not pointwise congruent.
Consequently, they constitute a Bonnet pair. Let us also remark that, by Proposition \ref{heli},
the surface $\psi^*$ is also properly helicoidal.

\end{proof}

It follows immediately from this proposition and the Lawson-type correspondence \eqref{fundadan} that there also exist Bonnet pairs in homogeneous $3$-manifolds with $\tau \neq 0$. More specifically, we have:

\begin{cor}
The sister surface in $\Ek$ of any simply connected properly helicoidal CMC surface in $\M^2 (\widetilde{\kappa})\times \R$ has a Bonnet mate in $\Ek$.
\end{cor}
This follows directly from the interesting fact that, by its own construction, the Lawson-type correspondence \eqref{fundadan} for CMC surfaces preserves positive Bonnet pairs, and by Proposition \ref{heli} it also preserves the property of being helicoidal.

\vspace{0.3cm}

{\bf Twin CMC immersions in $\Ek$ when $\tau \neq 0$:} Let $\psi :\Sigma\flecha \Ek$ be any simply
connected CMC surface in $\Ek$, $\tau \neq 0$, with fundamental data \eqref{fundamental}. Assume
also that $H\neq 0$. Consider now $\alfa\in \R$ so that $e^{i\alfa} (H+i\tau)= -H +i \tau$. It
follows then from the Lawson-type correspondence by Daniel exposed in formula \eqref{fundadan} that
there exists a surface $\psi^*:\Sigma\flecha \Ek$ with fundamental data $$ (\landa |dz|^2, u,
-H,e^{i\alfa} p\, dz^2, e^{i \alfa} A\, dz).$$ This surface, which by Theorem \ref{lasces} exists
and is unique up to congruences, was introduced in \cite{Dan}. Following that paper, we shall say
that $\psi^*$ is the \emph{twin immersion} of the CMC surface $\psi^*$, and also that $\psi$ and
$\psi^*$ are \emph{twin immersions}.

Our interest in twin immersions comes from the following fact, whose proof follows directly from Lemma \ref{datos}.

\begin{cor}
Any pair of twin CMC immersions in $\Ek$, $\tau \neq 0$, constitute a negative Bonnet pair in $\Ek$.
\end{cor}

\vspace{0.3cm}

{\bf Helicoidal mates in $\Ek$ when $\tau \neq 0$:} Similarly to Proposition \ref{helibon}, we have the following situation when $\tau \neq 0$:
 \begin{pro}\label{helimate}
Any non-minimal simply connected open piece of a helicoidal surface in $\Ek$, $\tau\neq 0$, has a negative Bonnet mate.
 \end{pro}
The proof of this result is similar to the one of Proposition \ref{helibon}, just taking into account that if an uple $(\landa,u,H,p,A)$ depends only on $s=(z+\bar{z})/2$ and verifies equations ${\bf (C.1)}$ to ${\bf (C.4)}$, then the uple $(\landa, u, -H,-\bar{p},-\bar{A})$ also has these properties. In addition, the surfaces are not pointwise congruent by Remark \ref{isome} and the condition $\tau\neq 0$.

\section{Proof of Theorem \ref{main1}}

Let $\psi,\psi^*:\Sigma\flecha \Ek$ be a Bonnet pair, and consider $\Sigma$ endowed with the
Riemann surface structure of the pair, as explained in the previous section. Then, by Lemma
\ref{datos} the fundamental data of $\psi$ and $\psi^*$ are of the form \eqref{fundamental} and \eqref{chulipiruli}, respectively.

Even though the present section is dedicated to the product spaces $\M^2 (\kappa)\times \R$, we will begin by considering the general situation of surfaces in $\Ek$, but with one restriction: \emph{we shall assume that $\ep =1$ if $\tau  \neq 0$,} i.e. that the Bonnet pair is positive.

We will work on a neighborhood $\Omega$ of a point $z\in \Sigma$ on which $u^* =\sigma u$ holds for a fixed $\sigma = \pm 1$. As the surfaces are real analytic, this condition thus holds everywhere automatically. We remark that in the smooth category, as $u^2 = (u^*)^2$ on our surface, the set of points $z_0$ with the property that $u^* =\sigma u$ holds around them is open and dense on $\Sigma$.

With all of this, it follows directly by ${\bf (C.2)}$ that  $(A-\sg A^*)_{\bar{z}} =0$, i.e. the $1$-form $\alfa dz = (A-\sg A^*)dz$ is holomorphic on $\Omega$.

Let us distinguish two cases. If $\alfa $ vanishes identically on $\Omega$,  then we have $u=\sg u^*$, $A^* =
\sg A$, and by ${\bf (C.1)}$ it also holds that $p^* = p$. By Remark \ref{isome} we see that $\psi|_{\Omega}$ and $\psi^*|_{\Omega}$ differ only by an isometry of $\Ek$. So $\psi|_{\Omega}$ and $\psi^*|_{\Omega}$ do not
constitute a Bonnet pair, and the same happens to $\psi$ and $\psi^*$ by real analyticity.

If $\alfa$ does not vanish identically on $\Omega$, then it has only  isolated zeros there. We shall work locally away from these zeros, as this will suffice for our purposes. 

With this, we can consider now around points in $\Omega$ with $\alfa \neq 0$ the new local conformal parameter $w= w(z) =i  \int^z \alfa$. By an abuse of notation, we shall also denote this new complex parameter by $z$. We then
have to bear in mind from now on that $\alfa =-i dz$ with respect to this parameter.

In addition, by ${\bf (C.4)}$ we see that $|A^*|^2 = |A|^2$. Thereby, there is a local real
function $f$ such that
 \begin{equation}\label{main1uno}
A= -f - \frac{i}{2} ,\hspace{1cm} \sg A^*= -f + \frac{i}{2}.
 \end{equation}
Hence, ${\bf (C.2)}$ is rewritten as
 \begin{equation}\label{main1dos}
-f_{\bar{z}} =\frac{\landa u}{2} (H+i\tau).
 \end{equation}
On the other hand, by differentiating ${\bf (C.3)}$ we see by means of ${\bf (C.1)}$, $u=\sg u^*$ and the fact that $\ep =1$ if $\tau \neq 0$, that

$$ \def\arraystretch{2} \begin{array}{lll}
 u_{z\bar{z}} & = & -H_{\bar{z}} A - (H-i \tau) A_{\bar{z}} - H_z \bar{A} - u (\kappa - 4 \tau^2) |A|^2 + \displaystyle \frac{2 p}{\landa^2} \landa_{\bar{z}} \bar{A} - \frac{2p}{\landa} \overline{A_z} \\ & = & \sg (- H_{\bar{z}} A^* - (H-i \tau) A^*_{\bar{z}} - H_z \bar{A^*} - \sigma u (\kappa - 4 \tau^2) |A^*|^2 + \displaystyle \frac{2 p^*}{\landa^2} \landa_{\bar{z}} \bar{A^*} - \frac{2p^*}{\landa} \overline{A^*_z}).
 \end{array} $$
Using in this last relation that $A- \sg A^*$ is holomorphic and the identity $|A|^2 = |A^*|^2$ we arrive at $$(A- \sg A^* ) H_{\bar{z}} + \overline{(A -\sg A^*)} H_z = -\frac{2 p}{\landa} (\overline{A_z} - \frac{\landa_{\bar{z}}}{\landa} \bar{A} ) + \frac{2\sigma \, p^*}{\landa} (\overline{A^*_z}  - \frac{\landa_{\bar{z}}}{\landa} \bar{A^*}).$$ At last, using here \eqref{main1uno} and ${\bf (C.0)}$ we get
 \begin{equation}\label{esta}
 i (H_z - H_{\bar{z}} )= \frac{2 u }{\landa } (|p^*|^2 - |p|^2).
 \end{equation}
Recalling now \eqref{dets}, and taking into account that $\psi,\psi^*$ have the same metric and principal curvatures, we conclude that $|p^*|^2 = |p|^2$. So, denoting $z=s+it$, \eqref{esta} turns into
 \begin{equation}\label{achese}
H_t =0, \hspace{1cm} \text{ that is, } \hspace{1cm} H = H(s).
 \end{equation}
In addition, by \eqref{main1dos} we know that $f_{\bar{z}} (H-i\tau)\in \R$. This is a first order linear PDE, and the fact that $H=H(s)$ lets us solve it explicitly. Indeed, we get
 \begin{equation}\label{efedelta}
 f=f(\delta), \hspace{1cm} \text{ where } \hspace{1cm} \delta = \tau \, t + \int H(s) ds.
 \end{equation}

From now on, we shall assume that $\tau =0$, and so $\Ek =\M^2 (\kappa)\times \R$. The case where $\tau \neq 0$ will be discussed in the next section.

Once here, suppose that $z_0$ has a neighbourhood $\mathcal{U}\subset \Omega$ on which $H$ vanishes identically. Then the surfaces $\psi|_{\mathcal{U}}, \psi^*|_{\mathcal{U}}$ are associate minimal surfaces in $\M^2 (\kappa)\times \R$, as explained in Section \ref{secmain}. Indeed, this follows directly from ${\bf (C.2)}$, ${\bf (C.3)}$ and ${\bf (C.4)}$. Thus, by real analyticity, both surfaces must globally be associated minimal surfaces.

Let us now work on an open piece of $\Omega$ on which $H$ does not vanish. As $\tau =0$, it follows from \eqref{efedelta} that $f=f(s)$, and therefore $A=A(s)$. By \eqref{main1dos} we also get $f' = -\landa u H$, where here $'$ denote differentiation with respect to the real parameter $s$. Consequently, $\landa u$ depends only on $s$. Thus, by ${\bf (C.4)}$ both $\landa$ and $u$ depend only on $s$, and finally, by ${\bf (C.3)}$, $p$ also depends only on $s$. In this argument we have used that $A$ cannot vanish identically on an open subset of a surface with non-vanishing $H$.

To sum up, we have obtained the following conclusion: around any point in $\Omega$ with $H\neq 0$ there exists a local conformal parameter $z$ such that the coefficients of the fundamental data of $\psi$ and $\psi^*$ depend only on $z+ \bar{z}$. By Proposition \ref{heli}, $\psi$ and $\psi^*$ are open pieces of helicoidal surfaces around such point. Moreover, as with respect to this specific local parameter we see that ${\rm Im} (A) \neq 0$ and ${\rm Im} (A^*) \neq 0$, we can infer that these helicoidal surfaces are in fact \emph{properly} helicoidal. As a matter of fact, we actually have $p^* = \bar{p}$ and $A^*= \bar{A}$. Again by real analyticity, the surfaces are, globally, open pieces of properly helicoidal surfaces. This completes the first part of the proof.

The converse has been already proved in Section 3, as we showed there the existence of associate minimal surfaces, as well as the existence of a Bonnet mate for any properly helicoidal surface in $\M^2 (\kappa)\times \R$ (Proposition \ref{helibon}).

\section{Proof of Theorem \ref{main2}}

We will divide the proof into two different cases, depending on whether the Bonnet pair is considered to be positive or negative.

\subsection*{Positive Bonnet pairs}

Let us start with a positive Bonnet pair, i.e. with the case where $\ep =1$, and keep the notations of the previous section. Thus, by \eqref{achese} and \eqref{efedelta} we have that $H=H(s)$ and $A=A(\delta)$. Also observe that as $\tau \neq 0$, the pair $(s,\delta)$ constitute local coordinates on the surface. Now, by \eqref{main1dos} and the relation $2 \delta_z = H-i\tau$ we get that $\landa u$ depend only on $\delta$. Consequently, by ${\bf (C.4)}$ we obtain that $\landa =\landa (\delta )$ and $u=u(\delta)$. Actually, 
 \begin{equation}\label{estre}
 \landa (\delta) u(\delta) = -f' (\delta).
 \end{equation}
Bearing this in mind, a direct computation shows that ${\bf (C.3)} $ translates into $$\frac{\landa (u' +2A)}{2\bar{A}} = \frac{-2p}{H-i\tau}.$$ This indicates that there is a smooth complex function $F=F(\delta)$ such that
 \begin{equation}\label{efegran}
 \frac{p}{H-i\tau} = F(\delta),  \hspace{1cm} F := -\frac{\landa (u'+ 2A)}{4\overline{A}}.
 \end{equation}
Let us consider two subcases: assume first that $H$ is constant on an open set (and thus everywhere, by real analyticity). Then, we see from the just discussed relations that all the fundamental data of $\psi$ and $\psi^*$ depend only on $\delta$ on this open set, and that $\delta = {\rm Re} (H-i\tau)(s+it)$ is a conformal parameter (since $H$ is constant). Thus, by Proposition \ref{heli} both $\psi$ and $\psi^*$ are CMC helicoidal surfaces in $\Ek$. Moreover, as the Lawson-type correspondence \eqref{fundadan} preserves helicoidal CMC positive Bonnet pairs, we conclude from Theorem \ref{main1} that these surfaces in $\Ek$ are exactly the sister CMC surfaces of the properly helicoidal CMC surfaces in $\M^2 (\widetilde{\kappa}) \times \R$, with $\widetilde{\kappa} =\kappa - 4\tau^2$.

Now assume that $H$ is not constant around some point. Then, by the Codazzi equation ${\bf (C.1)}$ and \eqref{efegran} we get the relation
 \begin{equation}\label{runi}
 F'(\delta) (H(s)^2 +\tau^2) = H'(s) (-F(\delta) +  \landa (\delta)/2) + (\kappa - 4\tau^2) \landa(\delta) u(\delta) A(\delta).
 \end{equation}
Let us examine this equation. First, observe that if we write $F(\delta )= F_1 (\delta) + i F_2 (\delta)$, then by decomposing \eqref{runi} into real and imaginary parts, we obtain using \eqref{estre} that 
 \begin{equation}\label{13cero}
 \left\{\def\arraystretch{1.5} \begin{array}{lll} (H(s)^2 +\tau^2) F_1' (\delta) & = & H' (s) (-F_1 (\delta) + \landa (\delta)/2) + (\kappa -4\tau^2) f(\delta) f' (\delta), \\ (H(s)^2 +\tau^2) F_2' (\delta) & = & -H' (s) F_2 (\delta) + (\kappa - 4\tau^2 ) f' (\delta) /2.
\end{array} \right.
 \end{equation}
If we work out the value of $f'(\delta)$ in the second equation, and afterwards substitute it into the first equation, we end up with 
 \begin{equation}\label{13uno}
(H(s)^2 +\tau^2) (F_1' (\delta) -2 f(\delta) F_2' (\delta)) = H' (s) (2f(\delta) F_2 (\delta) - F_1 (\delta) + \landa(\delta)/2). 
 \end{equation}
We are going to prove that under these circumstances, we must have that $F(\delta)$ is constant, and that both $\psi$ and $\psi^*$ are pointwise congruent.

In order to do so, assume first of all that $2f(\delta) F_2(\delta) - F_1(\delta) +\landa(\delta) /2 \neq 0$ on an open interval. It is then immediately inferred from \eqref{13uno} that 
 \begin{equation}\label{htg}
 H' (s)= a (H(s)^2 +\tau^2), \hspace{1cm} a\in \R.
 \end{equation}
If we now put together \eqref{runi} and \eqref{htg} we are left with $$ (H(s)^2 +\tau^2) (F' (\delta) +a F(\delta) -a \landa (\delta)/2)= (\kappa - 4\tau^2) f'(\delta) (f(\delta) + i/2).$$ So, as $H(s)$ is not constant, we see that $f'(\delta)=0$, i.e. that $A$ is constant. By ${\bf (C.2)}$ and ${\bf (C.4)}$ this implies that $u=0$ and $\landa$ is also constant. As a conclusion, $F(\delta)$ is also constant by \eqref{efegran}. 

On the other hand, suppose now that 
 \begin{equation}\label{13dos}
2f(\delta) F_2(\delta) - F_1(\delta) +\landa(\delta) /2 =0
 \end{equation}
along an interval, what implies by \eqref{13cero} that 
 \begin{equation}\label{13tres}
 F_1' (\delta) = 2f(\delta) F_2' (\delta)
 \end{equation}
If $F_2'(\delta)=0$, then $F(\delta)$ is constant. Otherwise, if we differentiate the second equation in \eqref{13uno} with respect to $s$, and separate the variables $s,\delta$ in the resulting equation, we conclude the existence of a real constant $a\neq 0$ such that 
 \begin{equation}\label{13five}
 2H(s) H'(s)= -a H''(s) \hspace{1cm} \text{ and } \hspace{1cm} F_2 (\delta)= a F_2 '(\delta).
 \end{equation}
Thus $H(s)^2 + a H'(s)= b$ for some $b\in \R$, and if we use \eqref{13five} in the second equation of \eqref{13uno} we can conclude the existence of constants $c,d\in \R$ with $c=2(\tau^2 +b)/(\kappa - 4\tau^2)$ such that 
 \begin{equation}\label{13six}
f(\delta)= c\, F_2 (\delta) +d.  \end{equation}
Now, by differentiating \eqref{13dos} and using \eqref{13tres} we see that $4f'(\delta) F_2 (\delta) + \landa' (\delta)=0$, what lets us conclude by means of \eqref{13five} and \eqref{13six} that 
 \begin{equation}\label{13siete}
 \landa (\delta)= -2 c F_2 (\delta)^2 + 2e, \hspace{1cm} e\in \R.
 \end{equation}
Next, observe that if we put together ${\bf (C.4)}$ with \eqref{main1uno} and \eqref{estre} we get the relation 
 \begin{equation}\label{13ocho}
 \landa(\delta) (4 f(\delta)^2 +1 )= \landa(\delta)^2 - f' (\delta)^2.
 \end{equation}
If in this relation we use \eqref{13five}, \eqref{13six} and \eqref{13siete} we conclude that $F_2 (\delta)$ is a root of a polynomial with constant coefficients of degree four. So, either $F_2 (\delta)$ (and hence $F(\delta)$) is constant, or all the coefficients of this polynomial must vanish. If this is the case, we would have the relations $$ - 2 c^3 = c^2, \hspace{.4cm} c^2 d =0,  \hspace{.4cm} 2c ( 4d^2 +1) - 8 e c^2 = 8 e c + \frac{c^2}{a^2}, \hspace{.4cm} cde =0,\hspace{.4cm} e (4d^2 +1)= 2e^2.$$ Firstly, if $c=0$, then by \eqref{13six} $f(\delta)$ is constant (i.e. $A(\delta)$ is constant), what proves by ${\bf (C.2)}$ that $u=0$, and consequently that $F(\delta)$ is also constant. On the other hand, if $c\neq 0$ we can conclude that $c= -1/2$, $d=0$, $e\in\{0,1/2\}$ and $2e-1 = 1/(4a^2)$. But the last relation indicates that $2e-1$ should be positive, and this is not possible.

So, we have concluded that $F(\delta)$ must be constant in order to have a Bonnet pair. Let us finally show that, even in this case, a Bonnet pair cannot occur, since the surfaces turn out to be pointwise congruent on an open set, and thus globally by real analyticity.

The condition that $F(\delta)$ is constant lets us translate \eqref{13cero} into
 \begin{equation}\label{13zero}
 \left\{\def\arraystretch{1.5} \begin{array}{lll} 0& = & H' (s) (-F_1 + \landa (\delta)/2) + (\kappa -4\tau^2) f(\delta) f' (\delta), \\ 0 & = & -H' (s) F_2 + (\kappa - 4\tau^2 ) f' (\delta) /2.
\end{array} \right.
 \end{equation}
It is deduced from the second equation that $f'(\delta)$ is constant. Taking this into account, we infer from the first equation in \eqref{13zero} that either $\landa = 2 F_1$ (and thus $f$ is constant), or else $f(\delta)/(-F_1 + \landa (\delta) /2) $ is constant. In any of the two cases, we get the existence of $a_1,a_2\in \R$ with 
 \begin{equation}\label{13B}
 f(\delta)= a_1 \landa (\delta) + a_2.
 \end{equation}
Using now \eqref{13B} and the constancy of $f' (\delta)$,  \eqref{13ocho} tells that $\landa (\delta)$ is the root of a third degree polynomial with constant coefficients and principal term $4a_1^2$. This shows by means of \eqref{13B} that $f(\delta)$ and $\landa (\delta)$ are necessarily constant. By ${\bf (C.2)}$ this proves that $u=0$, hence $\landa = 1+ 4f^2$ by ${\bf (C.4)}$. Now we can use \eqref{efegran} to infer that $F=-2A^2 = -2 (f+i/2)^2$. And at last, using that $F$ is constant in \eqref{efegran}, and that $H'(s)\neq 0$, the equation ${\bf (C.1)}$ simplifies to $2F= 1+ 4f^2$. This tells that $f=0$, from where $\sigma A^* = -A$ and thus $p(s) = p^*(s)$ if we use \eqref{efegran} for $\psi^*$ rather than for $\psi$. So, we conclude by Remark \ref{isome} that both surfaces are pointwise congruent on an open subset. Thus, they do not constitute a Bonnet pair.

To sum up, we have proved that a real analytic positive Bonnet pair in $\Ek$ for $\tau \neq 0$ is necessarily made up by two helicoidal CMC surfaces that are the sister immersions of two properly helicoidal CMC surfaces in some product space.

\subsection*{Negative Bonnet pairs}

Let us consider next the case of negative Bonnet pairs in $\Ek$, that is, the case $\ep = -1$. As was explained in Remark \ref{lodeachecero}, we can assume that $H$ never vanishes identically on an open set of $\Sigma$.

So, let $(\landa,u,H,p,A)$ and $(\landa, u^*, -H, p^*,A^*)$ be the coefficients of the fundamental data of a Bonnet pair $\psi,\psi^*:\Sigma\flecha \Ek$. Again, we will assume without loss of generality that $u^* =\sigma u$ holds for some $\sigma = \pm 1$.

From ${\bf (C.2)}$ we have $$ A_{\bar{z}} =\frac{\landa u}{2} (H-i\tau), \hspace{1cm} \sigma A^*_{\bar{z}} = \frac{\landa u}{2} (-H-i\tau),$$ and so $\overline{A_{\bar{z}}} + \sigma A^*_{\bar{z}} =0$. This implies the local existence of a complex function $\phi$ with $\phi_z =A$ and $\phi_{\bar{z}} = -\sigma \overline{A^*}$. We then have

\begin{equation}\label{nu1}
A-\sigma A^* =(\phi + \bar{\phi})_z :=\beta_z, \hspace{1.5cm} \beta := \phi + \bar{\phi}.
\end{equation}
If $A=\sigma A^*$ on some open set, then $u=0$ and, by ${\bf (C.2)}$, $A\, dz$ is a holomorphic $1$-form. So, we can choose a new complex parameter (that will still be denoted by $z$) so that $A=1$ holds on this open set (and thus $\landa =4$ by ${\bf (C.4)}$). By ${\bf (C.3)}$ we obtain that $p^* = - \bar{p}$, and at last ${\bf (C.1)}$ tells us that both $H$ and $p$ depend only on $z+\bar{z}$. Thus both surfaces are helicoidal. By composing one of the surfaces with an isometry $\Psi$ of $\Ek$ if necessary we get by Remark \ref{isome} that the surfaces are helicoidal mates. 

Let us work from now on in the open set $\mathcal{D} = \{z_0\in \Sigma : \beta_z (z_0)\neq 0\}$. By \eqref{nu1} and $|A|^2 = |A^*|^2$ we get

\begin{equation}\label{nu2}
A=\left(\frac{1}{2} + i f\right) \beta_z, \hspace{1.5cm} \sigma A^*=\left(-\frac{1}{2} + i f\right) \beta_z
\end{equation}
for some smooth real function $f$. Using \eqref{nu2} we can write ${\bf (C.2)}$ for $\psi$ and $\psi^*$ as the system

\begin{equation*}
\left\{\def\arraystretch{2} \begin{array}{lll}
if_{\bar{z}} \beta_z + \displaystyle\left(\frac{1}{2} + i f\right) \beta_{z\bar{z}}  & = & \displaystyle\frac{\landa u}{2} (H+i\tau), \\
if_{\bar{z}} \beta_z + \displaystyle\left(-\frac{1}{2} + i f\right) \beta_{z\bar{z}}  & = & \displaystyle\frac{\landa u}{2} (-H+i\tau),
\end{array}\right.
\end{equation*}
or equivalently,
\begin{equation}\label{nu3}
\left\{\def\arraystretch{1.5} \begin{array}{lll}
 \beta_{z\bar{z}}  & = & \landa H u, \\
f_{\bar{z}} \beta_z + f \beta_{z\bar{z}}  & = & \displaystyle\frac{\landa u}{2} \tau.
\end{array}\right.
\end{equation}
The second formula in \eqref{nu3} proves that $f_{\bar{z}} \beta_z \in \R$, i.e. $df \wedge d\beta =0$. As $\beta_z \neq 0$ on $\mathcal{D}$, we can infer that, locally,
 \begin{equation}\label{nu4}
 f=f(\beta).
 \end{equation}
We remark that if $f=0$ on some open set of $\Sigma$, then composing (if necessary) one of the surfaces with an isometry $\Psi$ of $\Ek$ as above, we get $A =\sigma A^*$, and this case has already been discussed. So, we will impose from now on the condition $f\neq 0$. 

If we differentiate ${\bf (C.3)}$, then using ${\bf (C.2)}$,  ${\bf (C.1)}$ and ${\bf (C.0)}$ we arrive at $$u_{z\bar{z}} = -H_{\bar{z}} A - (H^2 +\tau^2) \frac{\landa u}{2} - H_z \bar{A} - u |A|^2 (\kappa - 4\tau^2) - \frac{2|p|^2}{\landa} u.$$ If we now use that this expression must also hold for the fundamental data of $\psi^*$, as well as the relations $|A|=|A^*|$ and $|p|=|p|^*$, we get $$-H_{\bar{z}} A -H_z \bar{A} = \sigma (H_{\bar{z}} A^* + H_z \overline{A^*}),$$ which by \eqref{nu2} and $f\neq 0$ lets us conclude that $H_{\bar{z}}\beta_z = H_z\beta_{\bar{z}}$, or in other words,
 \begin{equation}\label{nu7}
 H=H(\beta) \hspace{1.5cm} \text{ locally. }
 \end{equation}
Moreover, from \eqref{nu3} and \eqref{nu4} we have
 \begin{equation}\label{nu8}
 \left(f(\beta) -\frac{\tau}{2 H(\beta)}\right) \beta_{z\bar{z}} + f'(\beta)|\beta_z|^2 =0.
 \end{equation}
We obtain then two different local possibilities:

\vspace{0.1cm}

{\bf Case A:} If $2 f(\beta) = \tau /H(\beta)$ on an interval, then \eqref{nu8} tells that $f$ (and thus $H$) is constant on an open piece of $\Sigma$, and $f=\tau/(2H)$. By \eqref{nu2} we get the existence of some $\alfa\in \R$ with $\sigma A^*= e^{i\alfa} A$. Moreover, again by \eqref{nu2} we can observe that $$e^{i\alfa} =\frac{-H +i\tau}{H+i\tau}.$$  Finally, if we write ${\bf (C.3)}$ for both $\psi$ and $\psi^*$, and we use that $A\neq 0$, we see that $p^* = e^{i\alfa} p$. Therefore, up to an isometry of $\Ek$, $\psi$ and $\psi^*$ are twin CMC immersions over an open subset of $\Sigma$ (and hence globally).

\vspace{0.1cm}

{\bf Case B:} If $2 f(\beta) - \tau /H(\beta)$ does not vanish along an interval, then \eqref{nu8} can be rephrased as follows: \emph{the function } $$ \mu (z)= F(\beta (z)), \hspace{1.5cm} F'(x)=\exp \left(\int \frac{f'(x)}{f(x) - (\tau/2 H(x))} dx \right) $$ \emph{is harmonic} on an open set of $\Sigma$. Thus $\mu$ has a local harmonic conjugate $\mu^*$ so that $\mu + i \mu^*$ is a new local conformal parameter on the surface. Let us keep denoting by $z=s+it$ this new parameter, but also taking into account that we now have 
 \begin{equation}\label{cagl}
  s=F(\beta) \hspace{1cm} \text{ and } \hspace{1cm}  F'(\beta) \beta_z=1.
 \end{equation}
Thus, in this new parameter, and using \eqref{nu4} and \eqref{cagl}, we obtain that \eqref{nu2} can be rewritten as
 \begin{equation}\label{nu11}
 A=g_1(s) + i g_2 (s), \hspace{1cm} \sigma A^* = -g_1(s) +i g_2(s) = -\bar{A},
 \end{equation}
for a pair of smooth real function $g_1,g_2$. By \eqref{nu7}, \eqref{cagl}, ${\bf (C.2)}$ and ${\bf (C.4)}$ we see that, locally,
 \begin{equation}\label{nu12}
 H=H(s),\hspace{1cm} \landa = \landa (s),\hspace{1cm} u=u(s).
 \end{equation}
Finally, using \eqref{nu11} and \eqref{nu12} in ${\bf (C.3)}$ as well as the fact that $u_z = u_{\bar{z}}$, we can conclude that $p=p(s)$ and $p^* = -\bar{p}$. Therefore $\psi$ and $\psi^*$ constitute a Bonnet pair, and both of them are (up to an isometry of $\Ek$) helicoidal mates in $\Ek$, as we described in Section \ref{secmain}

This finishes the uniqueness part of the proof. Existence was already shown in Section \ref{secmain}.

So, we have finished the proof of Theorem \ref{main2}.

\section{Final remarks}\label{last}

Let us start this section by inspecting when are two Bonnet mates in $\Ek$ globally congruent. For this, let us consider $\psi :\Sigma\flecha \Ek$ to be a helicoidal surface (with complete orbits) with fundamental data \eqref{fundamental}. By Proposition \ref{heli} there exists a conformal parameter $z=s+it$ for the surface such that $s$ varies in an interval $I\subset \R$ , that $t$ varies in all $\R$ and that the coefficients $(\landa,u,H,p,A)$ depend only on $s=(z+\bar{z})/2$.

Now, observe that if we reverse the orientation of $\psi$, its unit normal $\eta$ is replaced by $-\eta$, the parameter $s+it$ is replaced by $s-it$, and its fundamental data turn into $$(\landa |dz|^2, -u, -H, -\bar{p} \, dz^2, \bar{A} \, dz).$$ 

If $\tau =0$, and we consider $\Psi$ to be the symmetry in $\M^2 (\kappa)\times \R$ with respect to a vertical plane, it follows that $\psi^* (z):= \Psi (\psi (\bar{z}))$ has the fundamental data given by $(\landa |dz|^2 , u,H, \bar{p} \, dz^2, \bar{A}\, dz ).$ That is, $\psi^*$ is the Bonnet mate of the (properly) helicoidal surface $\psi$, and it is globally congruent to $\psi$.

If $\tau \neq 0$ and $\Psi$ is an isometry that reverses the orientations of both base and fiber, it can be shown in the same way that $\psi$ and $\psi^* (z):= \Psi (\psi (\bar{z}))$ constitute two globally congruent Bonnet mates.

In other words, we have proved: \emph{if $\cS_1$, $\cS_2$ are two helicoidal mates in $\Ek$, and $\cS_1$ is an open piece of the helicoidal surface $\overline{\cS_1}$ (with complete orbits), then $\cS_2$ is
congruent to another open piece of $\overline{\cS_1}$.} 

Nevertheless, it must also be emphasized that in general two Bonnet mates are not globally congruent. For instance, two associate minimal surfaces in $\M^2 (\kappa)\times \R$ related by an angle $\theta \neq 0,\pi$ are generically non-congruent. In the same way, the twin immersion of a simply connected CMC surface in $\Ek$, $\tau \neq 0$, is generically non-congruent to the original surface. For instance, if the fundamental data of the original CMC surface have no symmetries, its twin immersion cannot be globally congruent to it. We also point out that some helicoidal CMC examples are globally (and not pointwise) congruent to their twin immersions (see \cite{Dan}).

Apart from the above issue, let us consider briefly the compactness of Bonnet pairs. By Theorem \ref{main1}, it is immediate that \emph{a real-analytic compact surface in $\M^2(\kappa)\times \R$ cannot have a Bonnet mate}. This happens because there are no compact minimal surfaces in $\M^2(\kappa)\times \R$ other than the horizontal slices $h=h_0$ in $\S^2 (\kappa)$. Thus, any real-analytic compact surface in $\M^2 (\kappa)\times \R$ is always pointwise determined by its metric and principal curvatures.

This is no longer true for surfaces in $\Ek$, $\tau \neq 0$. This follows directly from Theorem \ref{main2} and the existence of (real-analytic) CMC spheres in these homogeneous spaces. It also follows from the fact that a rotational sphere in $\Ek$, $\tau \neq 0$ always has a helicoidal mate. Nevertheless, in these cases the Bonnet mates are globally congruent, by our discussion at the beginning of the section. Let us also remark that an important open problem in the classical $\R^3$ setting is whether there exist compact Bonnet pairs in $\R^3$.

We will finish by exposing some open problems related to the results of the present paper.

\begin{enumerate}
\item
Rigidity of CMC surfaces: As explained in the introduction of this paper, the classical Bonnet problem in $\R^3$ has two equivalent formulations, that are no longer equivalent when the target $3$-space has not constant curvature. So, it is natural to consider how many surfaces in $\Ek$ with the same metric and mean curvature function can exist. This problem seems to be quite involved, so it may be more interesting to consider its restriction to CMC surfaces: \emph{given a simply connected CMC-$H$ surface $S$ in $\Ek$, how many CMC-$H$ surfaces in $\Ek$ exist that are locally isometric to $S$?} This problem has been treated for minimal surfaces in $\M^2(\kappa)\times \R$ in \cite{HST}, but it is still unanswered even in that particular case.
 \item
Compact Bonnet pairs: it remains unsolved whether a real-analytic non-simply connected compact surface in $\Ek$ can have a Bonnet mate not globally congruent to it. By Theorem \ref{main1} and Theorem \ref{main2}, we know that if such surface exists, it must be a compact CMC surface, and its Bonnet mate must be its twin immersion. So, the problem reduces to the following question: \emph{does it exist a pair of twin compact CMC immersions with non-trivial topology in $\Ek$, $\tau \neq 0$?}. 
 \item
Other target $3$-manifolds: taking into account our results, it is natural to consider if the metric and principal curvatures are enough to determine pointwise a surface in a general Riemannian $3$-manifold. More specifically, there are two possible choices of ambient $3$-spaces for this generalized Bonnet problem that may be of special interest. One is the Riemannian $3$-space ${\rm Sol}_3$, i.e. the only homogeneous manifold belonging to the Thurston geometries that is missing in our discussion (since it has a $3$-dimensional isometry group). The other choice is the class $M^2 \times \R$ of Riemannian products of an abstract Riemannian surface with the real line, as this type of spaces generalize the homogeneous manifolds $\M^2 (\kappa)\times \R$, and are currently becoming a fashion research topic.
\end{enumerate}

\def\refname{References}

\end{document}